\documentclass[12pt]{article}
\usepackage{graphicx,amssymb}
\newcommand{\NOT}[1]{}

\newcommand{\figg}[3]{\begin{figure}[h]\centering\includegraphics[width=#2\textwidth]{fig#1.eps}\caption{#3}\label{fig:fig#1}\end{figure}}

\newtheorem{thmm}{Theorem}
\newtheorem{dfn}{Definition}

\title{Upper bounds for the piercing number of families of pairwise intersecting convex polygons}
\author{Meir Katchalski, Mathematics, Technion, Haifa\\ David Nashtir, Mathematics, Technion
Haifa}
\date{June 22, 2011}
\begin{document}
\maketitle
{\bf Abstract}. A convex polygon $A$ is {\em related} to a convex $m$-gon
$K= \bigcap_{i=1}^m k_i^+$, where $k_1^+,\ldots, k_m^+$ are the $m$ halfplanes
whose intersection is equal to $K$, if $A$ is the intersection of halfplanes
$a_1^+,\ldots,a_l$, each of which is a translate of one of the $k_i^+$-s.
The planar family ${\cal A}$ is {\em related} to $K$ if each $A \in {\cal A}$ is related to
$K$. We prove that any family of pairwise intersecting
convex sets related to a given $n$-gon has a finite piercing number
which depends on $n$.
In the general case we show $O(3^{n^3})$, while for a certain class of families, we decrease the bound to $4(n-2)$, and for $n=3,4$ the bound is 3 and 6 respectively.\\

{\dfn A convex polygon $P$ is {\em related} to a convex $m$-gon
$K= \bigcap_{i=1}^m k_i^+$, where $k_1^+,\ldots, k_m^+$ are the $m$ halfplanes
whose intersection is equal to $K$, if $P$ is the intersection of halfplanes
$a_1^+,\ldots,a_l^+$, each of which is a translate of one of the $k_i^+$-s. We use the convention that the line $l$ is the boundary of the halfplane $l^+$ and that $l^-$ is the
halfplane with boundary $l$ so that $l^+\cap l^-=l$.
The family ${\mathcal P}$ is {\em related} to $K$ if each $P \in {\mathcal P}$ is related to
$K$.}

{\thmm A convex family of pairwise intersecting sets related to an n-gon is $3^{n \choose 3}$ pierceable.}\\


{\thmm Let $\mathcal{F}$ be a family of pairwise intersecting sets
related to an n-gon $F$ with edges
$\overline{h},\overline{v},\overline{a}_1,\ldots,\overline{a}_{n-2}$,
such that $\overline{h}=[-1,0]$, $\overline{v}=[0,y]$ for any $y>0$,
and the edges $\overline{a}_1,\ldots,\overline{a}_{n-2}$ have
positive slopes. Then $\mathcal{F}$ is $4(n-2)$ pierceable.
If $n=3,4$ then the family is $3$ and $6$ pierceable respectively.}

{\dfn Let ${\mathcal P}$ be a family  related to $m$-gon $K= \bigcap_{i=1}^m k_i^+$.
A triangle $T$ is called empty or negative, if $T=\bigcap_{i=1}^3 l_i^-$, where  $l_1^-,l_2^-,l_3^-$ are minimal halfplanes that are translates of some $k_j^-$ ($j=1,\ldots,m$) such that $\bigcap_{i=1}^3 l_i^+= \emptyset$}.\\

{\bf Proof of theorem 1.} Let $\mathcal{F}$ be a family of pairwise intersecting
polygons related to a convex $n-$gon.  Observe the set of $n$ minimal halfplanes. Let $\mathcal{E}$
be the family of all empty triangles
created by them and let $N=e(n)=|\mathcal{E}|$. We will prove the theorem by induction on $N$. It's obviously true for $N=0$ since then the intersection of any three minimal halfplanes is not empty, hence the intersection of any three halfplanes is not empty,  hence by Helly's theorem  $\bigcap\mathcal{F}\neq \emptyset$. \\

Suppose $N>0$ and Let $E$ be an arbitrary triangle in  $\mathcal{E}$.
Observe the edges of $E$. Each of them comes from a line through an edge of some $F\in \mathcal{F}$.
Let $E={e}_1^-\cap {e}_2^- \cap {e}_{3}^-$, let $M_1,M_2,M_3$ be the midpoints of the edges of $E$, and let $M={m}_1^+\cap {m}_2^+ \cap {m}_{3}^+$ be the triangle
created by the midpoints, where $m_i$ is parallel to $e_i$ for $i=1,2,3$.\\

\NOT{
Now, let $\mathcal{F'} \subset \mathcal{F}$ be the subfamily of all polygons
that do not intersect with $\{M_1,M_2,M_3\}$, and for $i=1,2$ let
$F_i\in \mathcal{F'}$ be a polygon with the property:

\begin{equation} \label{property1} \mbox{ if } F_i \mbox{ has an edge } \overline{f}_i \mbox {
parallel to } m_i \mbox { then } e_i^+\subset {f_i}^+\subset m_i^+.
\end{equation}

Let ${\mathcal{F}}_i$ be the subfamily of all polygons with this
property. Since convexity and pairwise intersection imply that no
two lines through edges of any $F\in \mathcal{F}$ can intersect inside
$M$, it follows that for any $F_i\in {\mathcal{F}}_i$ and for $j=1,2,3;
j\neq i$ we have $m_j^+\subset {f_j}^+$. Hence, for any $F_i\in
{\mathcal{F}}_i$, we have $f_j^-\cap f_k^-\subset e_i^+\subset f_i^+$
which means that $f_j^-\cap f_k^-\cap f_i^-=\emptyset$. For any
$F_i\in \mathcal{F'}$ that does not have property~(\ref{property1}),
we immediately get $f_j^-\cap f_k^-\cap f_i^-=\emptyset$.\\

Thus, $\mathcal{F}$ splits into four subfamilies: one of them is pierced
by $\{M_1,M_2,M_3\}$ and the other three having at most $N-1$ empty
triangles each.

}

Since any two of sets intersect, for any $F \in {\mathcal{F}}$, $F=\bigcap_{i=1}^3 f_i^+$ contains at least one of the points $M_1, M_2, M_3$. Otherwise, for $i,j,k=1,2,3; i\neq j\neq k$, there would exist an angle $f_i^+\cap f_j^+$ which is strictly contained in the angle $m_i^+\cap m_j^+$, hence disjoint from the halfplane $e_k^+$, thus disjoint from a member of ${\mathcal{F}}$. It follows that  $\mathcal{F}$ can be divided into to three subfamilies, as follows:

$${\mathcal{F}}_1=\{F \in {\mathcal{F}} | M_1 \in \bigcap_{i=1}^3 f_i^+\}$$
$${\mathcal{F}}_2=\{F \in {\mathcal{F}} | M_2 \in \bigcap_{i=1}^3 f_i^+, M_1 \not\in \bigcap_{i=1}^3 f_i^+\}$$
$${\mathcal{F}}_3=\{F \in {\mathcal{F}} | M_3 \in \bigcap_{i=1}^3 f_i^+, M_1 \not\in \bigcap_{i=1}^3 f_i^+, M_2 \not\in \bigcap_{i=1}^3 f_i^+\}$$

where each of these subfamilies contains no empty triangle of type $E$, hence having at most $N-1$ empty triangles. Since by Helly's theorem, a family with no empty triangles is 1-pierceable, we get the following recursive inequality for the piercing number $f(N)$:


$$f(N) \leqslant 3f(N-1)$$

Hence:


$$f(N) \leqslant 3^N$$

Since the number of maximal empty triangles $N=e(n)< {n \choose 3}$
we finally get:


$$f(n) < 3^{n \choose 3}.\ \Box$$

{\bf More detailed explanations and drawings to be added later...}\\

{\bf Proof of theorem 2.} Let $\mathcal{F}$ be a family of sets related
to the convex $n-gon$ $F$ with edges
$\overline{h},\overline{v},\overline{a}_1,\ldots,\overline{a}_{n-2}$,
such that $\overline{h}=[-1,0]$, $\overline{v}=[0,y]$ for any $y>0$,
and the edges
$\overline{a}_1,\ldots,\overline{a}_{n-2}$ have positive slopes.\\

Observe the set of minimal halfplanes $h^+, v^+, a_i^+$ (
$i=1,\ldots,n-2$) and choose $1\leq s \leq n-2$), so that $\triangle E$
is the empty triangle $h^-\cap v^- \cap a_s^-\neq \emptyset$. Let $M_h, M_v$ and
$M_s$ be the midpoints of $\triangle E$, and let $\triangle M$ be
the triangle $\triangle M_h M_v M_s$ with edges $m_h, m_v,m_s$.
First, note that $\triangle M$
has the following:\\

{\bf Two Edges Outside (TEO) property.} Let  $F\in \mathcal{F}$ be a
polygon and let $L=\{l_1,l_2,l_3\}$ be a subset of its edges so that
$l_1\parallel h,l_2\parallel v,l_3\parallel a_s$. Then at most one
member of $L$ intersects $\triangle M$.\\

To prove the TEO property it is enough to notice that the pairwise
intersection implies that no polygon can have a vertex inside
$\triangle M$, hence no two edges of a polygon can intersect inside
$\triangle M$, hence if $l_i\in L$ intersects $\triangle M$, the
other two edges must lie outside $\triangle M$.\\

Now, proceed by choosing $a_s$ as the line with the smallest
positive slope (with respect to to the $x-axis$), such that 
$\triangle E=h^-\cap v^- \cap a_s^-\neq \emptyset$. 
Note that all polygons in $\mathcal{F}$ have edges $h',v'$ ($h'\parallel h,v'\parallel
v$), but there might exist ones that do not have edge $a'_s\parallel
a_s$. Let ${\mathcal{A}}_s\subset \mathcal{F}$ be the subfamily of all
polygons that have
edge $a'_s$ and do not intersect with $\{M_h, M_v, M_s\}$. We examine two cases.\\

{\bf Case 1.} $M_s\in a_i^+$ for $i=1,\ldots,n-2$.\\

We note that in this case, any polygon in ${\mathcal{A}}_s$ has the following properties:\\

\begin{equation} \label{property2}
a'_s \mbox{ is outside } \triangle M
\end{equation}

\begin{equation} \label{property3}
h' \mbox{ is outside } \triangle M
\end{equation}

\NOT{
\begin{equation} \label{property2}
\overline{M_vM_h}\in {a'_s}^+
\end{equation}

\begin{equation} \label{property3}
\overline{M_vM_s}\in {h'}^+
\end{equation}

first make a general observation
that because of pairwise intersection no two edges of a polygon can
intersect inside the midpoints triangle. Hence if a part of one of
the edges $h',v',a'$ lies inside the triangle, the other two must be
outside it. Thus,

}

To establish those properties, note that if property~\ref{property2}
does not hold, then by TEO both $h'$ and $v'$ are outside the triangle, hence,
both ${h'}^+$ and ${v'}^+$ contain $M_s$, and since  $M_s\in a_i^+$
for $i=1,\ldots,n-2$, it implies that the polygon itself
contains $M_s$ - a contradiction.\\

As for property~\ref{property3}, note that if  it does not hold,
then $h'$ intersects  $\triangle M$, and by TEO both ${v'}^+$ and ${a'}_s^+$ contain $M_h$. Further more, pairwise intersection implies that the intersection point $P=h'\cap v$ belongs
to the polygon, hence the intersection point $P_i=a'_i\cap h'$ for any $i=1,2,\ldots,n-2$ lies to the
left of $P$. Let $\alpha=\angle(a_s,h)$ and $\beta=\angle(a'_i,h)$  for $i\neq s$. 
If ${a'}_i^-$ does not create an empty triangle with $h^-$ and $v^-$ then ${a'}_i^+$ contains $M_h$. If, on the other hand, ${a'}_i^-\cap h^-\cap v^-\neq \emptyset$, then since $a_s$ has the smallest slope among all $a_i$'s that create an empty triangle with $h^-$ and $v^-$, it follows that $\beta>\alpha$, hence  again, $M_h\in {a'}_i^+$ - a contradiction. See figure~\ref{fig:fig9}.\\[10pt]

\figg{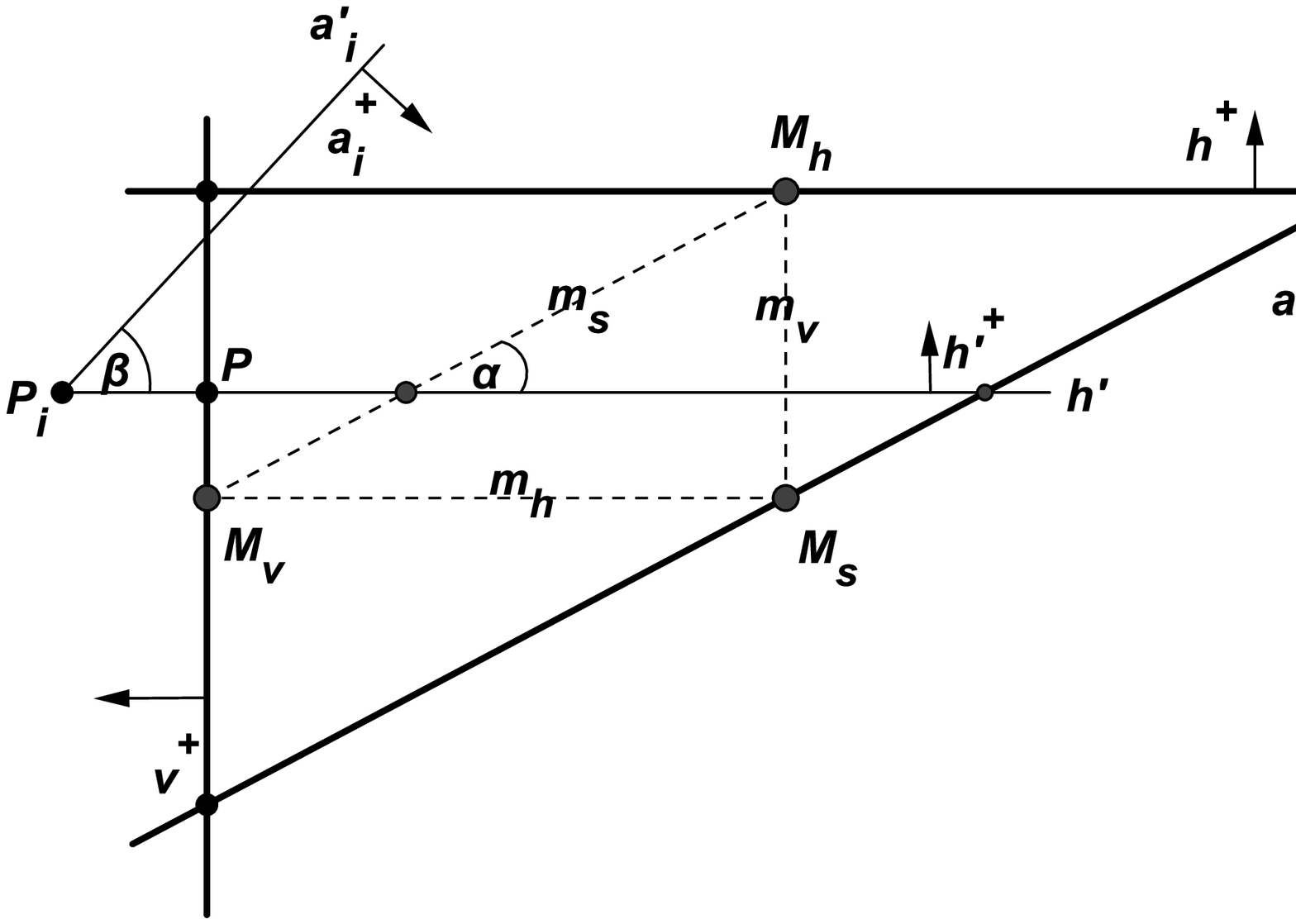}{1}{}

It follows that the members of ${\mathcal{A}}_s$ do not create empty triangles similar to $h^-\cap v^-\cap a_s^-$.\\

 {\bf Case 2.} There exits a
halfplane $a_i^-$ such that $M_s\in a_i^-$.\\

Let ${\mathcal{A}}^-(M_s)=\{a_i^-|M_s\in a_i^-\}$. Let $H=a_s\cap v$,
let $h_s$ be the horizontal line through $H$, let $v_s$ be the
vertical line through $M_s$ and let $P=v_s\cap h_s$. Choose an
arbitrary $a_i$ from ${\mathcal{A}}^-(M_s)$ and construct a new
auxiliary triangle $\triangle T_i=HXY$ as follows:\\

- if $P_i=a_i\cap h_s$ lies to the right of $P$, then $Y=P$ and
$X=M_s$.\\

- if $P_i=a_i\cap h_s$ lies to the left of $P$, then $Y=P_i$ and
$X=v_i\cap a_s$ where $v_i$ is the vertical line through $P_i$.\\

See figure~\ref{fig:fig8}.\\

\figg{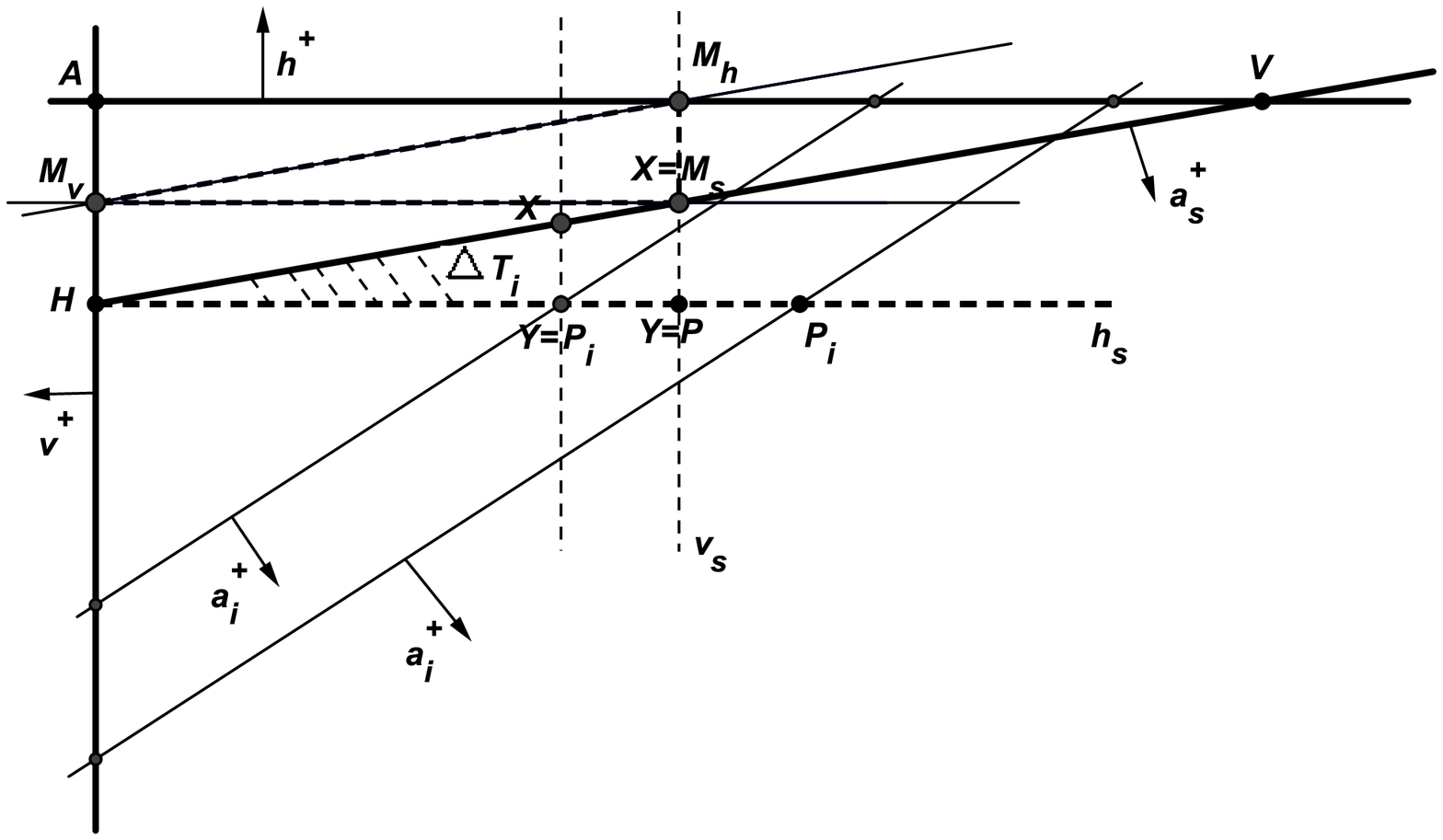}{1}{}

Denoting by
${\mathcal{A}}_i\subset \mathcal{F}$ the subfamily of all polygons that have
edge $a'_s$ and do not contain the vertex $X$, we see $\triangle T_i$ has properties similar to those of
$\triangle M_hM _vM_S$ we examined in case 1, i.e. any polygon in ${\mathcal{A}}_i$ has its $a'_s,h'$ edges outside $\triangle T_i$.\\

Hence if ${\mathcal{A}}_i\subset \mathcal{F}$ is the subfamily of all polygons that have
edge $a'_s$ and do not intersect with $\{M_h, M_v, M_s,X\}$ then the members of ${\mathcal{A}}$ do not create empty triangles similar to $h^-\cap v^-\cap a_s^-$.\\

Thus, we get the following recursive inequality for the
piercing number $f(N)$ where $N$ is the number of maximal empty triangles:

$$f(N) < f(N-1)+4$$

hence:

$$f(N) \leqslant 4N$$

and since $N\leqslant n-2$ we finally have $$f(n) \leqslant 4(n-2).\Box$$

\end{document}